\documentclass{amsart}
\usepackage{preamble}
\title{Integer decomposition property of three-dimensional smooth polytopes: a shorter proof}
\author{Juliana Curtis}
\author{Gaku Liu}

\begin{document}

\maketitle

\section{Main result}

We work in a real vector space $\bb R^d$ with fixed lattice $\bb Z^d$.
A lattice polytope $P \subset \bb R^d$ has the \emph{integer decomposition property} if for every positive integer $k$ and rational point $x \in P \cap \frac{1}{k} \bb Z^d$, there are lattice points $x_1$, \dots, $x_k \in P$ such that $x_1 + \dots + x_k = kx$.
A lattice polytope is \emph{smooth} if the direction cone at each of its vertices is generated by a lattice basis. Oda's smoothness conjecture \cite{oda2008problems} asserts that every smooth lattice polytope has the integer decomposition property.
We give a proof of the conjecture in dimension 3.

\begin{thm} \label{t:main}
Every three-dimensional smooth lattice polytope has the integer decomposition property.
\end{thm}

The first proof of this result was given in \cite{michalek26}, by AI. Here we give a shorter proof of the result. From the AI paper, we take one idea: the ``width lemma'' (Lemma 5.1 of \cite{michalek26}), a variation of which is our Lemma~\ref{width}. The remainder of the math, and all of the writing, was done without AI. We used Fable 5.1 to check the final draft for correctness and typos.

We would like to thank Mateusz Micha{\l}ek for providing us the result, and the positive nature of our discussions.

\begin{proof}

It is well-known \cite{ewald1991ampleness} that to prove the integer decomposition property for a given $d$-dimensional polytope, it suffices to prove the property for all $k \le d-1$. Thus, in $\bb R^3$ it suffices to prove the property for $k = 2$.

Let $P$ be a smooth lattice polytope in $\bb R^3$ and $x \in P \cap \frac{1}{2} \bb Z^3$. We want lattice points $x_1$, $x_2 \in P$ such that $x_1 + x_2  = 2x$. If $x$ is lattice height at most 2 from some facet of $P$, this follows from the next claim, which we prove in Section~\ref{sec:layer}.

\begin{lem} \label{layer}
If $x$ is lattice height at most 2 from some facet of $P$, then it is contained in a lattice polytope with the integer decomposition property contained in $P$.
\end{lem}

Thus, assume $x$ is lattice height at least $2.5$ from all facets of $P$. Let $R$ be the rational polytope
\[
R := P \cap (2x - P).
\]
The theorem is equivalent to saying that $R$ contains a lattice point. We will use the following ``width lemma'', proven in Section~\ref{sec:width}.

\begin{lem} \label{width}
$R$ has lattice width at least 4 with respect to every nonzero lattice functional.
\end{lem}

The result now follows by the bound in \cite{AVERKOV2021129} on the flatness constant in dimension 3, which implies that any convex body in $\bb R^3$ without a lattice point has width at most 3.972 with respect to some nonzero lattice functional. It remains to prove Lemmas~\ref{layer} and \ref{width}.
\end{proof}

\section{Proof of Lemma~\ref{layer}} \label{sec:layer}

Recall $P$ is a smooth lattice polytope in $\bb R^3$ and $x \in P \cap \frac{1}{2} \bb Z^3$ such that $x$ is lattice height at most 2 from some facet $F$ of $P$. If $x$ is lattice height at most 1 from $F$, then it is contained in the Cayley polytope $\conv(F,F')$ where $F'$ is the polygon obtained by shifting the supporting hyperplane of $P$ at $F$ inward by one lattice unit and intersecting with $P$. Since $P$ is smooth, $F'$ is a lattice polytope which is weakly Minkowski equivalent to $F$. By \cite{fakhruddin2002multiplicationmapslinearsystems} or \cite{Haase_2008}, this polytope has the integer decomposition property, so we are done.

Now assume $x$ is lattice height at least 1.5 from all facets of $P$. Let $Q$ be the polytope obtained by shifting all facets of $P$ one lattice unit inward. Then $x \in Q$. We prove the following.

\begin{lem} \label{contract}
$Q$ is a lattice polytope. The normal cone of each vertex of $Q$ satisfies one of the following.
\begin{enumerate}
    \item It is unimodular.
    \item It is lattice equivalent to the cone with generators $e_1$, $e_2$, $e_1 + e_2 - 2e_3$.
    \item It is lattice equivalent to the cone with generators $e_1$, $e_2$, $e_3$, $e_1+e_2-e_3$.
\end{enumerate}
\end{lem}

\begin{proof}
Since $P$ is smooth, the normal fan of $P$ has only these types of flips:
\begin{enumerate}[label=(\alph*)]
    \item $(2,2)$ flips which are lattice equivalent to
    \begin{multline*}
        \{e_1,e_2,e_3\}, \{e_1,e_2,me_1 + ne_2-e_3\} \to \\ \{e_1,e_3,me_1 + ne_2-e_3\}, \{e_2,e_3,me_1 + ne_2-e_3\}
    \end{multline*}
    for integers $m$, $n > 0$. The circuit involved is $(\{e_1,e_2\}, \{e_3,me_1 + ne_2-e_3\})$. Here and below, the circuit is ordered so that the flip removes the first part and adds the second part to the triangulation.
    \item $(1,3)$ flips which are lattice equivalent to
    \begin{multline*}
        \{e_1,e_2,e_3\}, \{e_1,e_3,-e_1-e_2+ne_3\}, \{e_2,e_3,-e_1-e_2+ne_3\} \to \\
        \{e_1,e_2,-e_1-e_2+ne_3\}
    \end{multline*}
    for some integer $n > 0$. The circuit involved is $(\{e_3\},\{e_1,e_2,-e_1-e_2+ne_3\})$.
    \item $(1,2)$ flips which are lattice equivalent to
    \begin{multline*}
        \{e_1,e_1+e_2,e_3\}, \{e_1,e_1+e_2,me_1 + ne_2-e_3\}, \\ \{e_2,e_1+e_2,e_3\}, \{e_2,e_1+e_2,me_1 + ne_2-e_3\} \to \\
        \{e_1,e_2,e_3\}, \{e_1,e_2,me_1 + ne_2-e_3\}
    \end{multline*}
    for some integers $m$, $n$. The circuit involved is $(\{e_1+e_2\}, \{e_1,e_2\})$.
    \item $(a,0)$ flips which decrease the dimension of $P$.
\end{enumerate}

The walls in the secondary fan of $P$ corresponding to these flips are given by the hyperplanes
\begin{equation} \label{circuit}
\sum_{i \in X_+} \lambda_i b_i = \sum_{i \in X_-} \lambda_i b_i,
\end{equation}
where $(X_+,X_-)$ is the circuit involved in the flip and the $\lambda_i$ are the absolute values of the coefficients in the linear relation of this circuit. (The $b_i$ are the coordinate system of the space in which the secondary fan lives, i.e., they are the right hand sides of the inequalities defining $P$.) If the $b_i$'s all decrease by the same amount, then the flip which replaces $X_+$ with $X_-$ can only occur if the right hand side of \eqref{circuit} decreases more than the left hand side, that is, $\sum_{i \in X_+} \lambda_i < \sum_{i \in X_-} \lambda_i$.

For any real $t \ge 0$, let $P_t$ be the polytope obtained by shifting each facet of $P$ inward by $t$ lattice units. As $t$ moves from 0 to 1, we track all possible flips from the list above.
\begin{enumerate}[label=(\alph*)]
    \item We have $\sum_{i \in X_+} \lambda_i = m+n$ and $\sum_{i \in X_-} \lambda_i = 2$. Thus this flip cannot occur as long as the normal fan of $P_t$ stays smooth because $m$, $n > 0$.
    \item We have $\sum_{i \in X_+} \lambda_i = n$ and $\sum_{i \in X_-} \lambda_i = 3$. This flip can only occur if $n < 3$. If $n = 2$, this flip can only occur at $t = 1$ because the value $\sum_{i \in X_+} \lambda_i b_i - \sum_{i \in X_-} \lambda_i b_i$ changes by 1 when $t$ changes by 1, and it is initially an integer. If $n = 1$, then the flip can occur at either $t = 1/2$ or $t = 1$, but the normal fan of $P_t$ remains smooth after the flip.
    \item We have $\sum_{i \in X_+} \lambda_i = 1$ and $\sum_{i \in X_-} \lambda_i = 2$. This flip can only occur at $t = 1$, by the same argument as above.
    \item These flips cannot occur because we know $x$ is lattice height at least $1.5$ from all facets, so $P_1 = Q$ is full-dimensional.
\end{enumerate}

Thus, the only flips which can change the smoothness of the normal fan of $P_t$ are flips (b) and (c) occurring exactly at $t = 1$. First, suppose that a flip occurs at $t=1$ and this flip does not share any full-dimensional simplices with another flip occurring. If the flip is of type (b), it will result in a cone of the form (1) or (2). If the flip is of type (c), it will produce two cones of type (1). Now suppose some flips occurring at $t=1$ share full-dimensional simplices. The only ways this can occur are
\begin{enumerate}[label=(\roman*)]
    \item A sequence of flips of type (c) whose circuits are in the same plane and are of the form $(\{v_2\}, \{v_1,v_3\})$, $(\{v_3\},\{v_2,v_4\})$, $(\{v_4\}, \{v_3,v_5\})$, \dots, for some vectors $v_1$, $v_2$, \dots. The first two circuits of this chain must be lattice equivalent to $(\{e_1+e_2\}, \{e_1,e_2\})$ and $(\{e_2\}, \{e_1+e_2,-e_1\})$. The sum of the linear relations of these two circuits is the linear relation $e_1 + (-e_1) = 0$. Hence, if these two flips occur simultaneously, then a flip of type (d) with circuit $(\{e_1,-e_1\},\emptyset)$ also occurs. This is a contradiction, so no such chain exists. 
    \item Two flips of type (c) with circuits lattice equivalent to $(\{e_1+e_2\}, \{e_1,e_2\})$ and $(\{e_1+e_2\}, \{e_3,e_1+e_2-e_3\})$. This results in a cone of type (3).
    \item Two flips of type (c) with circuits lattice equivalent to $(\{e_1+e_2\}, \{e_1,e_2\})$ and $(\{e_3\}, \{e_1+e_2,-e_1-e_2+e_3\})$. This results in two cones of type (1).
    \item A flip of type (b) and a flip of type (c) with circuits lattice equivalent to $(\{e_3\},\{e_1,e_2,-e_1-e_2+ne_3\})$ and $(\{-e_1-e_2+ne_3\}, \{e_3,-e_1-e_2+(n-1)e_3\})$, where $n=1$ or 2. If $n = 2$, then this results in a cone of type (1). If $n = 1$, then the sum of the linear relations of the two circuits gives the linear relation $e_1 + e_2 + (-e_1-e_2) = 0$. Hence, if these two flips occur simultaneously, then a flip of type (d) with circuit $(\{e_1,e_2,-e_1-e_2\},\emptyset)$ also occurs, a contradiction.
\end{enumerate}

In conclusion, all cones are of the claimed form. Since the normal fan of $P_t$ is smooth for all $t < 1$, the support function of $P_1 = Q$ is linear on a smooth fan. In addition, it takes integer values on the primitive rays of this fan. Thus $Q$ is a lattice polytope.
\end{proof}

We return to the proof of Lemma~\ref{layer}. We first claim $x$ is lattice height at most 1 from some facet $G$ of $Q$. Let $F$ be a facet of $P$ such that the lattice height from $x$ to $F$ is minimal. By assumption, this height is at most 2. If $Q$ has a facet with the same outer normal as $F$, then $x$ is lattice height at most 1 from this facet and we are done. Otherwise, the primitive outer normal to $F$ disappears from the normal fan when moving from $P$ to $Q$. This can only happen due to a $(1,b)$ flip ($b = 2$ or 3).  In this case, the disappearing ray $v$ can be written as $v = \sum_{i=1}^b c_i v_i$ where $v_i$ are different primitive rays of the normal fan, and $\sum c_i > 1$. Let $h_t$ be the support function of $P_t - x$, with $P_t$ defined as in the proof of Lemma~\ref{contract}. At the time $t$ of the aforementioned flip, we have $h_t(v) = \sum_{i=1}^b c_i h_t(v_i)$. Since $\sum c_i > 1$ and $h_t(v) > 0$, we must have $h_t(v_i) < h_t(v)$ for some $i$. Then $h_0(v_i) < h_0(v)$, contradicting our assumption that the lattice height from $x$ to $F$ is minimal.

Thus $x$ is lattice height at most 1 from some facet $G$ of $Q$. By Lemma~\ref{contract}, the direction cone of $Q$ at each vertex of $G$ is one of the following.
\begin{enumerate}
    \item A unimodular cone.
    \item Lattice equivalent to the cone with generators $2e_1+e_3$, $2e_2+e_3$, $e_3$.
    \item Lattice equivalent to the cone with generators $e_1$, $e_2$, $e_1+e_3$, $e_2+e_3$.
\end{enumerate}

Let $G'$ be the rational polygon obtained by shifting the supporting hyperplane of $Q$ at $G$ inward by one unit and intersecting with $Q$. Let $G''$ be the convex hull of all lattice points in $G'$. Every vertex $v$ of $G'$ is also a vertex of $G''$, unless $v$ is on an edge from a vertex $u$ of $G$ of type (2). In this case, without loss of generality let $u = 0$ and $G$ be contained in the plane $\spa(2e_1+e_3,2e_2+e_3)$.\footnote{We may assume this because the cone in (2) is symmetric with respect to $\GL(\bb Z^3)$.} Then $v = \frac{1}{2}e_3$, and $G''$ instead has vertices $e_1+e_3$ and $e_2+e_3$, which are points on the two edges of $G'$ from $v$.

Let $A = \conv(G,G')$ and $B = \conv(G,G'')$. We have $x \in A$. If $x \in A \setminus B$, then it is contained in some tetrahedron lattice equivalent to $\conv(0,2e_1+e_3,2e_2+e_3,e_3)$. This has the integer decomposition property because it is a tetrahedron with height 1 with respect to its facet $\conv(2e_1+e_3,2e_2+e_3,e_3)$.
So assume $x \in B$. We will use the following result, which follows from the proof of Theorem 3.2 in \cite{Haase_2008}.

\begin{thm}[\cite{Haase_2008}] \label{cayley}
Let $X$, $Y$ be lattice polygons with the following property: For every edge $e$ of $Y$ and $\ell$ a primitive lattice functional which achieves its maximum over $Y$ on $e$, there exists a lattice line segment $s$ in $X$ which is parallel to $e$, such that $\ell$ achieves either its maximum or one less than its maximum over $X$ on $s$.
Then the Cayley sum $\Cay(X,Y)$ has the integer decomposition property.
\end{thm}

Let $e$ be an edge of $G''$ and $\ell$ a primitive lattice functional which achieves its maximum over $G''$ on $e$. Then $\ell$ achieves its maximum over $G$ at either an edge of $G$ or a vertex $u$ of type (2) or (3). 
\begin{enumerate}[start=2]
    \item If $u$ is of type (2), without loss of generality let $u = 0$ and let $G$ be contained in the plane $\spa(2e_1+e_3,2e_2+e_3)$. As noted before, $e$ is the segment $[e_1 + e_3,e_2 + e_3]$. Then $G$ contains the parallel line segment $[2e_1+e_3,2e_2+e_3]$, and $\ell$ achieves one less than its maximum over $G$ on this segment.
    \item If $u$ is of type (3), without loss of generality let $u = 0$ and let $G$ be contained in the plane $\spa(e_1,e_2)$. Then $e = [e_1+e_3,e_2+e_3]$. Then $G$ contains the parallel line segment $[e_1,e_2]$, and $\ell$ achieves one less than its maximum over $G$ on this segment.
\end{enumerate}
Thus, the conditions of Theorem~\ref{cayley} are satisfied, so $\conv(G,G'')$ has the integer decomposition property, completing the proof.

\section{Proof of Lemma~\ref{width}} \label{sec:width}

The proof here is a strengthening of the argument in Lemma 5.1 of \cite{michalek26}.

Recall that $P$ is a smooth lattice polytope in $\bb R^3$, $x \in P \cap \frac{1}{2} \bb Z^3$ such that $x$ is lattice height at least $2.5$ from all facets of $P$, and $R = P \cap (2x-P)$. For any convex set $K$, let $h_K : \bb R^3 \to \bb R$ denote its support function. Let
$R' := R - x$, so $R'$ is symmetric about the origin.
We wish to prove that $h_{R'}(v) \ge 2$ for all $v \in \bb Z^3 \setminus \{0\}$.

Let $P' = P - x$. By our assumptions on $x$, we have $h_{P'}(v) \ge 2.5$ for all facet normals $v \in \bb Z^3$ of $P$. In fact, since $P$ is smooth, for all $v \in \bb Z^3$ we can write $h_{P'}(v) = \sum_{i=1}^3 c_i h_{P'}(v_i)$ for some facet normals $v_i$ and nonnegative integers $c_i$, and thus
\begin{equation} \label{firstbound}
h_{P'}(v) \ge 2.5 \text{ for all } v \in \bb Z^3 \setminus \{0\}.
\end{equation}

\begin{lem} \label{segment}
Let $u_0$ and $v$ be lattice vectors with $v \neq 0$, and $n$ a nonnegative integer such that the sequence $(h_{P'}(u_0 + iv))_{i=0}^n$ is a weakly increasing arithmetic sequence with difference $d$. Then $h_{P'}(u_0+nv) \ge 2.5n$, and if $d < 2.5$, then $h_{P'}(u_0) \ge 2.5n$ and $h_{P'}(u_0+nv) \ge (2.5+d)n$.
\end{lem}

\begin{proof}
If $n = 0$ or $d \ge 2.5$, then the result follows immediately since $h_{P'}(u) \ge 0$ for all $u \in \bb Z^3$ by \eqref{firstbound}. So assume $n \ge 1$ and $d < 2.5$.

We first claim $u_0 \neq 0$. Suppose $u_0 = 0$. Let $e_1$, $e_2$, $e_3$ be the primitive generators of a full-dimensional normal cone of $P$ containing $v$. Then $v = \sum_{i=1}^3 c_i e_i$ for nonnegative integers $c_i$ (since $P$ is smooth), and
\[
d = h_{P'}(v) = \sum_{i=1}^3 c_i h_{P'}(e_i) \ge 2.5
\]
by \eqref{firstbound}. This contradicts $d < 2.5$. Thus $u_0 \neq 0$.

If $n = 1$, then the result follows from $h_{P'}(u_0) \ge 2.5$. Assume $n \ge 2$. Set $u_i := u_0 + iv$. Since $h_{P'}$ is convex, and forms an arithmetic progression on $u_0$, \dots, $u_n$ where $u_1 \neq u_0$, $u_n$, it follows that $h_{P'}$ is affine on $[u_0,u_n]$. Hence $[u_0,u_n]$ is contained in some normal cone of $P$.

Let $e_1$, $e_2$, $e_3$ be the primitive generators of a full-dimensional normal cone containing $[u_0,u_n]$. Let $v = \sum_{i=1}^3 c_i e_i$, where the $c_i$ are integers, not all zero. Then $d = \sum_{i=1}^3 c_i h_{P'}(e_i)$. Since $d < 2.5$, there is some $c_i$ which is negative. Without loss of generality, assume $c_1 \le -1$. Write $u_0 = \sum_{i=1}^3 \alpha_i e_i$ and $u_n = \sum_{i=1}^3 \beta_i e_i$ where $\alpha_i$ and $\beta_i$ are nonnegative integers. Then $\alpha_1 = \beta_1 - nc_1 \ge n$. Hence $h_{P'}(u_0) \ge 2.5n$. The inequality $h_{P'}(u_n) \ge (2.5+d)n$ follows immediately.
\end{proof}

\begin{lem} \label{growth}
Let $u_0$ and $v$ be lattice vectors with $v \neq 0$, and $n$ a nonnegative integer such that the sequence $(h_{P'}(u_0 + iv))_{i=0}^n$ is weakly increasing. Then $h_{P'}(u_0 + nv) \ge 2n - 0.5n_0$, where $n_0$ is the largest $i$ such that $h_{P'}(u_0) = h_{P'}(u_0+iv)$.
\end{lem}

\begin{proof}
Set $u_i := u_0 + iv$. Let $d_i := h_{P'}(u_i) - h_{P'}(u_{i-1})$. We have $h_{P'}(u_i) = h_P(u_i) - \langle u_i,x \rangle$, and hence
\[
d_i = h_{P}(u_i) - h_{P}(u_{i-1}) - \langle v, x \rangle.
\]
Since $h_P(u_i)$ and $h_P(u_{i-1})$ are integers, it follows that either $d_i \in \bb Z$ for all $i$ or $d_i \in \bb Z + \frac{1}{2}$ for all $i$.

First assume $d_i \in \bb Z$ for all $i$. Since $h_{P'}$ is convex, the sequence $d_1$, \dots, $d_n$ is weakly increasing. Let $n_0$, $n_1$, $n_{\ge 2}$ be the numbers of indices $i$ for which $d_i = 0$, $d_i = 1$, and $d_i \ge 2$, respectively. We have the two separate inequalities
\begin{align*}
    h_{P'}(u_n) &\ge h_{P'}(u_{n_0}) + n_1 + 2n_{\ge 2} \\
    h_{P'}(u_n) &\ge h_{P'}(u_{n_0+n_1}) + 2n_{\ge 2}.
\end{align*}
Applying Lemma~\ref{segment} to the sequences $\{u_i\}_{i=0}^{n_0}$ and $\{u_i\}_{i=n_0}^{n_0+n_1}$, we obtain
\begin{align*}
    h_{P'}(u_n) &\ge 2.5n_0 + n_1 + 2n_{\ge 2} \\
    h_{P'}(u_n) &\ge 3.5n_1 + 2n_{\ge 2}.
\end{align*}
Adding $3/5$ the first inequality to $2/5$ the second inequality, we have
\[
h_{P'}(u_n) \ge 1.5n_0 + 2n_1 + 2n_{\ge 2} = 2n - 0.5n_0
\]
as desired.

Now assume $d_i \in \bb Z + \frac{1}{2}$. Let $n_{0.5}$, $n_{1.5}$, $n_{\ge 2.5}$ be the numbers of indices $i$ for which $d_i = 0.5$, $d_i = 1.5$, and $d_i \ge 2.5$, respectively. By a similar argument as above, we have
\begin{align*}
    h_{P'}(u_n) &\ge 3n_{0.5} + 1.5n_{1.5} + 2.5n_{\ge 2.5} \\
    h_{P'}(u_n) &\ge 4n_{1.5} + 2.5n_{\ge 2.5}.
\end{align*}
Adding $4/5$ the first inequality to $1/5$ the second inequality, we have
\[
h_{P'}(u_n) \ge 2.4n_{0.5} + 2n_{1.5} + 2.5n_{\ge 2.5} \ge 2n
\]
as desired.
\end{proof}

\begin{lem} \label{sum}
Let $u_0$ and $v$ be lattice vectors with $v \neq 0$, and $n$ a nonnegative integer. Then $h_{P'}(u_0) + h_{P'}(u_0+nv) \ge 2n$.
\end{lem}

\begin{proof}
Set $u_i := u_0 + iv$ and $d_i := h_{P'}(u_i) - h_{P'}(u_{i-1})$ as before. Let $a$ be the smallest index such that $d_{a+1} \ge 0$ (with $a = n$ if there is no such index), and let $b$ be the largest index such that $d_b \le 0$ (with $b = 0$ if there is no such index). Let $n_0 = b-a$ be the number of indices $i$ with $d_i = 0$. Applying Lemma~\ref{growth} to the sequences $\{u_i\}_{i=a}^n$ and $\{u_{b-i}\}_{i=0}^b$, we obtain
\[
    h_{P'}(u_n) + h_{P'}(u_0) \ge 2(n + n_0) - n_0 \ge 2n.
\]
\end{proof}

We can now complete the proof of Lemma~\ref{width}. Let $v \in \bb Z^3 \setminus \{0\}$. We have $R' = (-P') \cap P'$. By linear program duality, we have
\begin{align}
h_{R'}(v) &= \min_{\substack{w_0,w_1 \in \bb Q^3 \\ w_0 + w_1 = v}} ( h_{-P'}(w_0) + h_{P'}(w_1) ) \nonumber \\
&= \min_{\substack{w_0,w_1 \in \bb Q^3 \\ -w_0 + w_1 = v}} ( h_{P'}(w_0) + h_{P'}(w_1) ). \label{dual}
\end{align}
Let $w_0$, $w_1 \in \bb Q^3$ such that the minimum in \eqref{dual} is achieved.
Let $n$ be a positive integer such that $nw_0$, $nw_1 \in \bb Z^3$. Then $nw_1 = nw_0 + nv$, so by Lemma~\ref{sum}, we have $h_{P'}(nw_0) + h_{P'}(nw_1) \ge 2n$. Hence $h_{R'}(v) = \frac{1}{n}(h_{P'}(nw_0) + h_{P'}(nw_1)) \ge 2$, as desired.

\printbibliography

\end{document}